\documentclass{elsart}

\usepackage{ifpdf}
\usepackage{graphicx,natbib,amssymb,lineno, amsmath, mathrsfs}
\begin{document}

\begin{frontmatter}
\hsize=6.5in


\title{New identities for the Glasser transform and their applications}

\author[label1]{Ahmet Dernek},
\author[label1]{Ne\c se Dernek},
\author[label2]{Osman Y\"urekli\corauthref{cor1}}

\address[label1]{Department of Mathematics, University of Marmara, TR-34722, Kadik\"oy, Istanbul, Turkey}
\address[label2]{Department of Mathematics, Ithaca College, Ithaca, NY 14850, US}

\corauth[cor1]{Corresponding author}

\begin{abstract}

In the present paper the authors show that  an iteration of the $\mathscr{L}_{2}$-transform by itself is a constant multiple of the Glasser transform. Using this iteration identity, a Parseval-Goldstein type theorem for $\mathscr{L}_{2}$-transform and the Glasser transform is given. By making use of these results  a number of new Parseval-Goldstein type identities are obtained for these and many other well-known integral transforms. The identities proven in this paper are shown to give rise to useful corollaries for evaluating infinite integrals of special functions. Some examples are also given as illustration of the results presented here.

\end{abstract}


\end{frontmatter}

\section{Introduction}
\numberwithin{equation}{section}
Over a  decade ago, \citet{YS} presented a systematic account
of so-called the $\mathscr{L}_{2}$-transform:
\begin{equation}
\label{l2}
\mathscr{L}_{2}\big\{f(x)\, ;\, y\big\}
=\int_0^\infty x\,\exp\big(-x^{2}\,y^{2}\big)\,f(x)\,dx
\end{equation}
The $\mathscr{L}_{2}$-transform is related to the classical Laplace transform 
\begin{equation}
\label{lap}
\mathscr{L} \big\{f(x)\, ;\, y \big\}
=
\int_0^\infty  \exp(-x \,y )\,f(x)\,dx
\end{equation}
by means of the following relationships:
\begin{align}
\label{l2:l}
\mathscr{L}_{2}\left\{f(x)\, ;\, y\right\}
&=
\frac{1}{2}\,\mathscr{L}\Big\{f\big(\sqrt{x}\big)\, ;\, y^{2}\Big\},
\\
\label{ll2}
\mathscr{L}\left\{f(x)\, ;\, y\right\}
&=
2\,\mathscr{L}_{2}\Big\{f\big(x^{2}\big)\, ;\, \sqrt{y}\Big\}.
\end{align}
Subsequently, various Parseval-Goldstein type identities were given in (for example)  \citet{BDY:1}, \citet{BDY:2}, \citet{DSY}, \citet{Gl}, \citet{Y99a}, and \citet{Y99b} for the $\mathscr{L}_2$-transform. New solutions techniques were obtained for the Bessel differential equation in \citet{WY02} and the Hermite differential equation in \citet{WY03} using this integral transform. There are numerous analogous results in the literature on various integral transforms (see, for instance  \citet{Y89}, \citet{SY95},\citet{YG}, and \citet{YS98}). Some of the results from \citet{Y89} and \citet{Y92} are applied to generalized functions by \citet{AA}.

Over three decades ago, \citet{Gl} considered so-called the Glasser transform
\begin{equation}
\label{glasser}
\mathscr{G}\big\{f(x)\, ;\, y\big\}
=
\int_{0}^{\infty}
\frac{f(x)}{\sqrt{x^{2}+y^{2}}}\,dx.
\end{equation}
Glasser gave the following Parseval-Goldstein type theorem (cf. \citet[p. 171, Eq. (4)]{Gl})
\begin{equation}
\label{pggl}
\int_{0}^{\infty}f(x)\,\mathscr{G}\big\{g(y)\, ;\, x\big\}\,dx
=
\int_{0}^{\infty}g(x)\,\mathscr{G}\big\{f(y)\, ;\, x\big\}\,dx,
\end{equation}
and evaluated a number of infinite integrals involving Bessel functions. Additional results about the Glasser transform
can be found in \citet{SY95} and  \citet{KSY}.

The Fourier sine transform and the Fourier cosine transform are defined as 
\begin{align}
\label{fsin}
\mathscr{F}_{S}\big\{f(x) \, ;\, y \big\}&=\int_0^\infty \sin(x\,y)\,f(x)\,dx,
\intertext{and}
\label{fcos}
\mathscr{F}_{C}\big\{f(x) \, ;\, y \big\}&=\int_0^\infty \cos(x\,y)\,f(x)\,dx,
\end{align}
respectively.

  The Hankel transform is defined by
\begin{equation}
\label{hankel}
\mathscr{H}_{\nu}\big\{f(x)\, ;\, y\big\}
=
\int_{0}^{\infty}
\sqrt{x\,y}\,\text{J}_{\nu}(x\,y)\,f(x)\,dx
\end{equation}
where ${\rm J}_{\nu}(x)$ is the Bessel function of the first kind of order $\nu$. Using the formula (cf. \citet[p. 306, Eq. 32:13:10]{SO})
\begin{equation}
\label{bessel:half}
{\rm J}_{1/2}(x)=\sqrt{\frac{2}{\pi\,x}}\,\sin(x),
\end{equation}
the definition \eqref{fsin} of the Fourier sine transform, and the definition \eqref{hankel} of the Hankel transform, we obtain the familiar relationship 
\begin{equation}
\label{hankel:sin}
\mathscr{H}_{1/2}\big\{f(x) \, ;\, y \big\}=\sqrt{\frac{2}{\pi}}\,\mathscr{F}_{S}\big\{f(x) \, ;\, y \big\}\cdot
\end{equation}
Similarly, using the formula (cf. \citet[p. 306, Eq. 32:13:11]{SO})
\begin{equation}
\label{bessel:nhalf}
\text{J}_{-1/2}(x)=\sqrt{\frac{2}{\pi\,x}}\,\cos(x),
\end{equation}
the definition \eqref{fcos} of the Fourier cosine transform, and the definition \eqref{hankel} of the Hankel transform, we obtain the relationship 
\begin{equation}
\label{hankel:cos}
\mathscr{H}_{-1/2}\big\{f(x) \, ;\, y \big\}=\sqrt{\frac{2}{\pi}}\,\mathscr{F}_{C}\big\{f(x) \, ;\, y \big\}\cdot
\end{equation}
The $\mathscr{K}$-transform is defined by
\begin{equation}
\label{kn}
\mathscr{K}_{\nu}\big\{f(x)\, ;\, y\big\}
=
\int_{0}^{\infty}
\sqrt{x\,y}\,{\rm K}_{\nu}(x\,y)\,f(x)\,dx
\end{equation}
where ${\rm K}_{\nu}$ is the Bessel function of the second kind of order $\nu$. Using the formula (cf. \citet[p. 239, Eq. 26:13:5]{SO})
\begin{equation}
\label{kbessel:half}
{\rm K}_{1/2}(x)=\sqrt{\frac{\pi}{2\,x}}\exp(-x),
\end{equation}
the definition \eqref{lap} of the Laplace transform, and the definition \eqref{kn} of the $\mathscr{K}$-transform, we obtain the relationship 
\begin{equation}
\label{kn:e}
\mathscr{K}_{1/2}\big\{f(x) \, ;\, y \big\}=\sqrt{\frac{\pi}{2}}\,\mathscr{L}\big\{f(x) \, ;\, y \big\},
\end{equation}
which incidentally holds true also when $\mathscr{K}_{1/2}$ is replaced by $\mathscr{K}_{-1/2}$.

In this article, we show that an iteration of the $\mathscr{L}_{2}$-transform by itself is a constant multiple of the Glasser transform defined by \eqref{glasser}. Using this iteration identity, we establish a Parseval-Goldstein type theorem relating the $\mathscr{L}_{2}$-transform and the Glasser transform. The Parseval-Goldstein type theorem established here yields potentially new identities for the various integral transform introduced above. As applications of the resulting identities and theorems, some illustrative examples are also given. 

\section{The Main Theorem}
\numberwithin{equation}{section}
In the following lemma, we give an iteration identity involving the $\mathscr{L}_{2}$-transform \eqref{l2} and the Glasser transform \eqref{glasser}.
\begin{lem}\label{l1} 
The identity 
\begin{equation}
\label{l1:1}
\mathscr{L}_{2}\Big\{\frac{1}{u}\,\mathscr{L}_{2}\Big\{f(x)\, ;\, u\Big\}\, ;\, y\Big\}
=
\frac{\sqrt{\pi}}{2}
\,\mathscr{G}\big\{x\,f(x)\, ;\, y\big\},
\end{equation}
holds true, provided that the integrals involved converge absolutely.
\end{lem}
\begin{pf} 
Using the definition \eqref{l2} of the $\mathscr{L}_{2}$-transform,
we have
\begin{equation}
\label{l1:p1} 
\mathscr{L}_{2}\Big\{\frac{1}{u}\,\mathscr{L}_{2}\big\{f(x)\, ;\, u\big\}\, ;\, y\Big\}
=
\int_{0}^{\infty} 
\exp\big(-y^{2}\,u^{2}\big)\,
\bigg[
\int_{0}^{\infty} 
x\,\exp\big(-x^{2}\,u^{2}\big)\,f(x)\,dx
\bigg]\,du.
\end{equation}
Changing the order of integration, which is permissible by absolute convergence of the integrals involved, and then using the definition \eqref{l2} of the $\mathscr{L}_{2}$-transform once more, we find from \eqref{l1:p1} that
\begin{align}
\notag
\mathscr{L}_{2}\Big\{\frac{1}{u}
\,\mathscr{L}_{2}\big\{f(x)\, ;\, u\big\}\, ;\, y\Big\}
&=
\int_{0}^{\infty} 
x\,f(x)\,
\bigg[
\int_{0}^{\infty} 
\exp\Big[\big(-y^{2}+x^{2}\big)u^{2}\Big]\,du\bigg]\,dx
\\
\label{l1:p2}
&=
\int_{0}^{\infty} 
x\,f(x)\,
\mathscr{L}_{2}\Big\{\frac{1}{u}\,\, ;\, \big(x^{2}+y^{2}\big)^{1/2}\Big\}
\,dx.
\end{align}
Furthermore, we have
\begin{equation}
\label{l1:p3}
\mathscr{L}_{2}\Big\{\frac{1}{u}\,\, ;\, \big(x^{2}+y^{2}\big)^{1/2}\Big\}=\frac{\sqrt{\pi}}{2}\,\big(x^{2}+y^{2}\big)^{-1/2}.
\end{equation}
Now the assertion \eqref{l1:1} follows from \eqref{l1:p2}, \eqref{l1:p3}, and the definition \eqref{glasser} of the Glasser-transform.
\qed
\end{pf}
The Lemma \ref{l1} yields some useful corollaries that will be required in our investigation.
\begin{cor}\label{c1:l1}
We have (cf. \citet[p. 171, (2)]{Gl})
\begin{equation}
\label{c1:l1:1}
\mathscr{G}\big\{x^{\mu-1} \, ;\, y \big\}=2^{-\mu}\,{\rm B}\Big(\mu,\frac{1}{2}-\frac{\mu}{2}\Big)\,y^{\mu-1}, \quad 0<\Re(\mu)<1,
\end{equation}
where ${\rm B}(x,y)$ is the beta function defined by
\begin{equation}
\label{beta}
{\rm B}(x,y)=\int_{0}^{1}t^{x-1}(1-t)^{y-1}\,dt, \quad x>0, y>0,
\end{equation}
and it is related to the gamma function through
\begin{equation}
\label{b:g}
{\rm B}(x,y)=\frac{\Gamma(x)\,\Gamma(y)}{\Gamma(x+y)}={\rm B}(y,x).
\end{equation}
\end{cor}
\begin{pf}
We set
\begin{equation}
\label{c1:l1:p1}
f(x)=x^{\mu-2}, \quad 0<\Re(\mu)<1
\end{equation}
in Lemma \ref{l1}. Using the relation \eqref{l2:l} and the known formula \citet[p. 137, Entry (1)]{E1}, we find that
\begin{equation}
\label{c1:l1:p2}
\mathscr{L}_{2}\big\{x^{\mu-2} \, ;\, u \big\}=\frac{1}{2}\,\mathscr{L}\big\{x^{(\mu-2)/2} \, ;\, u^{2} \big\}
=\frac{1}{2}\,\Gamma\Big(\frac{\mu}{2}\Big)\,u^{-\mu}.
\end{equation}
Multiplying the equation \eqref{c1:l1:p2} through by $1/u$ and than applying the $\mathscr{L}_{2}$-transform, we obtain
\begin{equation}
\label{c1:l1:p3}
\mathscr{L}_{2}\Big\{\frac{1}{u}\,\mathscr{L}_{2}\big\{x^{\mu-2} \, ;\, u \big\}\, ;\, y\Big\}=
\frac{1}{2}\,\Gamma\Big(\frac{\mu}{2}\Big)\,\mathscr{L}_{2}\big\{u^{-\mu-1},y\big\}\cdot
\end{equation}
Using the relation \eqref{l2:l} and the formula \citet[p. 137, Entry (1)]{E1} once more on the right hand side of \eqref{c1:l1:p3}, we deduce that 
\begin{equation}
\label{c1:l1:p4}
\mathscr{L}_{2}\Big\{\frac{1}{u}\,\mathscr{L}_{2}\big\{x^{\mu-2} \, ;\, u \big\}\, ;\, y\Big\}=
\frac{1}{4}\,\Gamma\Big(\frac{\mu}{2}\Big)\,\Gamma\Big(\frac{1}{2}-\frac{\mu}{2}\Big)\,y^{\mu-1}.
\end{equation}
Utilizing the well-known duplication formula for the gamma function (cf. \citet[p. 414, Eq. (43:5:7)]{SO})
\begin{equation}
\label{c1:t1:p5}
\Gamma(2\alpha)=\frac{4^{\alpha}}{2\sqrt{\pi}}\,\Gamma(\alpha)\,\Gamma\Big(\frac{1}{2}+\alpha\Big)
\end{equation}
with $\mu=2\alpha$, the relationship \eqref{b:g} for the beta function on the right hand side of \eqref{c1:l1:p4} and finally the identity \eqref{l1:1} of our Lemma \ref{l1}, we obtain  the desired result \eqref{c1:l1:1}.
\qed
\end{pf}
\begin{cor}\label{c2:l1}
We have (cf. \citet[p. 171, (h)]{Gl})
\begin{equation}
\label{c2:l1:1}
\mathscr{G}\big\{x^{\nu+1}\,{\rm J}_{\nu }(z\,x) \, ;\, y \big\}
=\sqrt{\frac{2}{\pi\,z}}\,y^{\nu+\frac{1}{2}}\,{\rm K}_{\nu+\frac{1}{2}}(z\,y).
\end{equation}
\end{cor}
\begin{pf}
We set 
\begin{equation}
\label{c2:l1:p1}
f(x)=x^{\nu}\,{\rm J}_{\nu}(z\,x),\quad -1<\Re(\nu)<\frac{1}{2}
\end{equation}
in Lemma 1. Using the relation \eqref{l2:l} and the known formula \citet[p. 185, Entry (30)]{E1}, we find that
\begin{align}
\notag
\mathscr{L}_{2}\big\{x^{\nu}\,{\rm J}_{\nu}(z\,x)\, ;\, u\big\}
&=\frac{1}{2}
\mathscr{L}\Big\{x^{\nu/2}\,{\rm J}_{\nu}\big(z\,x^{1/2}\big)\, ;\, u^{2}\Big\}\\
\label{c2:l1:p2}
&=\frac{1}{2}\,\Big(\frac{z}{2}\Big)^{\nu}\,u^{-2\nu-2}\,\exp\bigg(-\frac{z^{2}}{4u^{2}}\bigg)\cdot
\end{align}
Multipliying the equation \eqref{c2:l1:p2} through by $1/u$ and than applying the $\mathscr{L}_{2}$-transform, we obtain
\begin{equation}
\label{c2:l1:p3}
\mathscr{L}_{2}\Big\{\frac{1}{u}\,\mathscr{L}_{2}\big\{x^{\nu}\,{\rm J}_{\nu}(z\,x)\, ;\, u\big\} \, ;\, y \Big\}
=\frac{1}{2}\,\Big(\frac{z}{2}\Big)^{\nu}\,
\mathscr{L}_{2}\bigg\{u^{-2\nu-3}\,\exp\bigg(-\frac{z^{2}}{4u^{2}}\bigg)\, ;\, y\bigg\}
\end{equation}
We evaulate the $\mathscr{L}_{2}$-transform on the right hand side of \eqref{c2:l1:p3} by using  the relation \eqref{l2:l} and the known formula \citet[p. 146, Entry (29)]{E1}:
\begin{align}
\notag
\mathscr{L}_{2}\bigg\{u^{-2\nu-3}\,\exp\bigg(-\frac{z^{2}}{4u^{2}}\bigg)\, ;\, y\bigg\}
&=\frac{1}{2}\,\mathscr{L}\bigg\{u^{-(2\nu+3)/2}\,\exp\bigg(-\frac{z^{2}}{4u}\bigg)\, ;\, y^{2}\bigg\}
\\\label{c2:l1:p4}
&=
\Big(\frac{2y}{z}\Big)^{\nu+\frac{1}{2}}\,{\rm K}_{\nu+\frac{1}{2}}(z\,y).
\end{align} 
Now the assertion \eqref{c2:l1:1} immediately follows upon substituting the result \eqref{c2:l1:p4} into the equation \eqref{c2:l1:p3} and using the identity \eqref{l1:1} of our Lemma \ref{l1}.
\qed
\end{pf}
\begin{cor}\label{c3:l1}
We have (cf. \citet[p. 174, (g)]{Gl})
\begin{equation}
\label{c3:l1:1}
\mathscr{G}\big\{{\rm J}_{\nu }(z\,x) \, ;\, y \big\}
={\rm I}_{\nu/2}\Big(\frac{1}{2}\,z\,y\Big)\,{\rm K}_{\nu/2}\Big(\frac{1}{2}\,z\,y\Big), \quad\Re(\nu)>-1.
\end{equation}
\end{cor}
\begin{pf}
We set 
\begin{equation}
\label{c3:l1:p1}
f(x)=\frac{1}{x}\,{\rm J}_{\nu}(z\,x),\quad \Re(\nu)>-1
\end{equation}
in Lemma \ref{l1}. Using the relation \eqref{l2:l} and the known formula \citet[p. 185, Entry 29]{E1}, we find that
\begin{equation}
\label{c3:l1:p2}
\mathscr{L}_{2}\Big\{\frac{1}{x}\,{\rm J}_{\nu}(z\,x)\, ;\, u\Big\}=\frac{\sqrt{\pi}}{2u}
\exp\bigg(-\frac{z^{2}}{8u^{2}}\bigg)\,{\rm I}_{\nu/2}\bigg(\frac{z^{2}}{8u^{2}}\bigg)\cdot
\end{equation}
Multipliying the equation \eqref{c3:l1:p2} through by $1/u$ and than applying the $\mathscr{L}_{2}$-transform, we obtain
\begin{equation}
\label{c3:l1:p3}
\mathscr{L}_{2}\bigg\{\frac{1}{u}\,\mathcal{L}_{2}\Big\{\frac{1}{x}\,{\rm J}_{\nu}(z\,x)\, ;\, u\Big\} \, ;\, y \bigg\}
=\frac{\sqrt{\pi}}{2}\,
\mathscr{L}_{2}\bigg\{\frac{1}{u^{2}}\,\exp\bigg(-\frac{z^{2}}{8u^{2}}\bigg)\,{\rm I}_{\nu/2}\bigg(\frac{z^{2}}{8u^{2}}\bigg)\, ;\, y\bigg\}\cdot
\end{equation}
Now the assertion \eqref{c3:l1:1}  immediately follows upon using the relation \eqref{l2:l} once more and than utilizing the known formula \citet[p. 325, Entry 10]{PBM} and \eqref{l1:1} of our Lemma \ref{l1}.
\qed
\end{pf}
\begin{thm}\label{t1} 
If the conditions stated in Lemma \ref{l1} are satisfied, then the Parseval-Goldstein type relations
\begin{align}
\label{t1:1}
\int_{0}^{\infty}
\mathscr{L}_{2}\left\{f(x)\, ;\, y\right\}\,\mathscr{L}_{2}\left\{g(u)\, ;\, y\right\}\,dy
&=
\frac{\sqrt{\pi}}{2}
\int_{0}^{\infty}
x\,f(x)\,
\mathscr{G}\big\{u\,g(u)\, ;\, x\big\}\,dx
\\
\label{t1:2}
\int_{0}^{\infty}
\mathscr{L}_{2}\left\{f(x)\, ;\, y\right\}\,\mathscr{L}_{2}\left\{g(u)\, ;\, y\right\}\,dy
&=
\frac{\sqrt{\pi}}{2}
\int_{0}^{\infty}
u\,g(u)\,
\mathscr{G}\big\{x\,f(x)\, ;\, u\big\}\,du
\\
\intertext{and}
\label{t1:3}
\int_{0}^{\infty}
x\,f(x)\,
\mathscr{G}\big\{u\,g(u)\, ;\, x\big\}\,dx
&=
\int_{0}^{\infty}
u\,g(u)\,
\mathscr{G}\big\{x\,f(x)\, ;\, u\big\}\,du
\end{align} 
hold true.
\end{thm}
\begin{pf} 
We only give the proof of \eqref{t1:1}, as the proof of \eqref{t1:2} is similar. Identity \eqref{t1:3} follows from the identities \eqref{t1:1} and \eqref{t1:2}. 

Using the definition \eqref{l2} of the $\mathscr{L}_{2}$-transform, we have
\begin{align}
\notag
\int_{0}^{\infty}
\mathscr{L}_{2}&\big\{f(x)\, ;\, y\big\}\,\mathscr{L}_{2}\big\{g(u)\, ;\, y\big\}\,dy
\\
\label{t1:p1}
&=
\int_{0}^{\infty}
\mathscr{L}_{2}\big\{g(u)\, ;\, y \big\}\,
\left[\int_{0}^{\infty}x\,\exp\big(-x^{2}\,y^{2}\big)\,f(x)\,dx\right]\,dy.
\end{align} 
Changing the order of integration (which is permissible by absolute convergence of the integrals involved) and using the definition \eqref{l2} of the $\mathscr{L}_{2}$-transform once again, we find from \eqref{t1:p1} that
\begin{align}
\notag
\int_{0}^{\infty}
\mathscr{L}_{2}&\big\{f(x)\, ;\, y\big\}\,\mathscr{L}_{2}\big\{g(u)\, ;\, y\big\}\,dy
\\
\notag
&=
\int_{0}^{\infty}
x\,f(x)
\left[\int_{0}^{\infty}\exp\big(-x^{2}\,y^{2}\big)\,
\mathscr{L}_{2}\big\{g(u)\, ;\, y\big\}\,dy\right]\,dx
\\
\label{t1:p2}
&=
\int_{0}^{\infty}
x\,f(x)\,
\mathscr{L}_{2}\bigg\{\frac{1}{y}\,
\mathscr{L}_{2}\big\{g(u)\, ;\, y\big\}\, ;\, x \bigg\}\,dx.
\end{align} 
Now the assertion \eqref{t1:1} easily follows from \eqref{t1:p2} and \eqref{l1:1} of the Lemma~\ref{l1}.
\qed
\end{pf}
\begin{cor}\label{c1:t1}
If the integrals involved converge absolutely and $0<\Re(\mu)<1$, then we have
\begin{align}
\label{c1:t1:1}
\int_{0}^{\infty} y^{-\mu}\,\mathscr{L}_{2}\big\{f(x)\, ;\, y\big\}\,dy
&=\frac{1}{2}\,\Gamma\Big(\frac{1}{2}-\frac{\mu}{2}\Big)
\,\int_{0}^{\infty}x^{\mu}\,f(x)\,dx,\\
\label{c1:t1:2}
\int_{0}^{\infty} y^{-\mu}\,\mathscr{L}_{2}\big\{f(x)\, ;\, y\big\}\,dy
&=\frac{\sqrt{\pi}}{\Gamma(\mu/2)}
\int_{0}^{\infty}u^{\mu-1}\,\mathscr{G}\big\{x\,f(x)\, ;\, u\big\}\,du, 
\intertext{and}
\label{c1:t1:3}
\int_{0}^{\infty}u^{\mu-1}\,\mathscr{G}\big\{x\,f(x)\, ;\, u\big\}\,du
&=
\frac{1}{2}\,{\rm B}\Big(\frac{\mu}{2},-\frac{\mu}{2}+\frac{1}{2}\Big)
\,\int_{0}^{\infty}x^{\mu}\,f(x)\,dx.
\end{align}
\end{cor}
\begin{pf}
We start with the proof of the assertion \eqref{c1:t1:1} by setting 
\begin{equation}
\label{c1:t1:p1}
g(u)=u^{\mu-2}
\end{equation}
in the Theorem \ref{t1}. Utilizing the formulas \eqref{c1:l1:1}, \eqref{c1:l1:p2} and the  identity \eqref{t1:1} of the Theorem \ref{t1}, we obtain that
\begin{equation}
\label{c1:t1:p2}
\int_{0}^{\infty}y^{-\mu}\,\mathscr{L}_{2}\big\{f(x)\, ;\, y\big\}\,dy
=\frac{\sqrt{\pi}}{2^{\mu}}\,\bigg[\Gamma\Big(\frac{\mu}{2}\Big)\bigg]^{-1}\,
B\Big(\mu,\frac{1}{2}-\frac{\mu}{2}\Big)\,\int_{0}^{\infty}x^{\mu}\,f(x)\,dx.
\end{equation}
Using the  duplication formula \eqref{c1:t1:p5} for the gamma function
with $\mu=2\alpha$ on the right hand side of \eqref{c1:t1:p2} we deduce  the assertion \eqref{c1:t1:1}.

Similarly, the proof of the assertion \eqref{c1:t1:2} follows upon utilizing \eqref{c1:t1:p1} and 
\eqref{c1:l1:p2} into the identity \eqref{t1:2} of our Theorem \ref{t1}. 

Finally, the assertion \eqref{c1:t1:3} easily follows from the identities \eqref{c1:t1:1}, \eqref{c1:t1:2} and the relationship \eqref{b:g} between the beta function and the gamma function. 
\qed
\end{pf}
\begin{cor}\label{c2:t1} 
If the integrals involved converge absolutely and $-1<\Re(\nu)<1/2$, then we have
\begin{align}
\label{c2:t1:1}
\mathscr{L}_{2}
\bigg\{y^{2\nu-1}\,\mathscr{L}_{2}\Big\{f(x)\, ;\, \frac{1}{2y}\Big\}\, ;\, z\bigg\}
&=2^{-\nu-\frac{1}{2}}\,z^{-\nu-1}\,\mathscr{K}_{\nu+\frac{1}{2}}\big\{x^{\nu+1}\,f(x)\, ;\, z\big\},
\\
\label{c2:t1:2}
\mathscr{L}_{2}
\bigg\{y^{2\nu-1}\,\mathscr{L}_{2}\Big\{f(x)\, ;\, \frac{1}{2y}\Big\}\, ;\, z\bigg\}
&=\frac{\sqrt{\pi}}{2^{\nu+1}}\,z^{-\nu-\frac{1}{2}}\,\mathscr{H}_{\nu}\Big\{u^{\nu+\frac{1}{2}}\,\mathscr{G}\big\{x\,f(x)\, ;\, u\big\}\, ;\, z\Big\},
\intertext{and}
\label{c2:t1:3}
\mathscr{K}_{\nu+\frac{1}{2}}\big\{x^{\nu+1}\,f(x)\, ;\, z\big\}&=
\Big(\frac{\pi\,z}{2}\Big)^{1/2}\,\mathscr{H}_{\nu}\Big\{u^{\nu+\frac{1}{2}}\,\mathscr{G}\big\{x\,f(x)\, ;\, u\big\}\, ;\, z\Big\},
\end{align}
where $\mathscr{H}_{\nu}\big\{f(x)\, ;\, y\big\}$ and $\mathscr{K}_{\nu}\big\{f(x)\, ;\, y\big\}$ denote the Hankel transform and the $\mathscr{K}$-transform as defined by \eqref{hankel} and \eqref{kn}, respectively.
\end{cor}
\begin{pf}
We put
\begin{equation}
\label{c2:t1:p1}
g(u)=u^{\nu}\,J_{\nu}(z\,u)
\end{equation}
in our Theorem \ref{t1}. Utilizing the identity \eqref{c2:l1:1} of Corollary \ref{c2:l1}, the equation \eqref{c2:l1:p2} and the Parseval-Goldstein type relation \eqref{t1:1} of Theorem \ref{t1}, we obtain
\begin{equation}
\label{c2:t1:p2}
\int_{0}^{\infty}\frac{1}{y^{2\nu+2}}\,\exp\bigg(-\frac{z^{2}}{4y^{2}}\bigg)\,\mathscr{L}_{2}\big\{f(x)\, ;\, y\big\}\,dy=\Big(\frac{2}{z}\Big)^{\nu+\frac{1}{2}}\,\int_0^\infty x^{\nu+\frac{3}{2}}\,{\rm K}_{\nu+\frac{1}{2}}(z\,x)\,f(x)\,dx.
\end{equation}
Making a simple change of variable in the integral on the left hand side and using the definition \eqref{kn} of the $\mathscr{K}$-transform on the right hand side of \eqref{c2:t1:p2} we obtain the desired identity \eqref{c2:t1:1}.

The assertion \eqref{c2:t1:2} is obtained similarly using the Parseval-Goldstein relation \eqref{t1:2} of Theorem \ref{t1}. The assertion \eqref{c2:t1:3} immediately follows from the relations \eqref{t1:1} and \eqref{t1:2}.
\qed
\end{pf}
\begin{rem}\label{r1:t1} 
If we let $\nu=0$ in our Corollary \ref{c2:t1} and than use the formula \eqref{bessel:half} and the definition \eqref{lap} of the Laplace transform, we obtain
\begin{align}
\label{r1:t1:1}
\mathscr{L}_{2}
\bigg\{\frac{1}{y}\,\mathscr{L}_{2}\Big\{f(x)\, ;\, \frac{1}{2y}\Big\}\, ;\, z\bigg\}
&=\frac{\sqrt{\pi}}{2z}\,\mathscr{L}\big\{x\,f(x)\, ;\, z\big\}\\
\label{r1:t1:2}
\mathscr{L}_{2}
\bigg\{\frac{1}{y}\,\mathscr{L}_{2}\Big\{f(x)\, ;\, \frac{1}{2y}\Big\}\, ;\, z\bigg\}
&=\frac{1}{2}\,\sqrt{\frac{\pi}{z}}\,\mathscr{H}_{0}\Big\{\sqrt{u}\,\mathscr{G}\big\{x\,f(x)\, ;\, u\big\}\, ;\, z\Big\}
\intertext{and}
\label{r1:t1:3}
\mathscr{L}\big\{x\,f(x)\, ;\, z\big\}
&=
\sqrt{z}\,\mathscr{H}_{0}\Big\{\sqrt{u}\,\mathscr{G}\big\{x\,f(x)\, ;\, u\big\}\, ;\, z\Big\}\cdot
\end{align}
\end{rem}
\begin{rem}\label{r2:t1} 
If we let $\nu=-1/2$ in our Corollary \ref{c2:t1} and than use the formula \eqref{hankel:cos} and the definition \eqref{fcos} of the Fourier cosine transform, we obtain
\begin{align}
\label{r2:t1:1}
\mathscr{L}_{2}
\bigg\{\frac{1}{y^{2}}\,\mathscr{L}_{2}\Big\{f(x)\, ;\, \frac{1}{2y}\Big\}\, ;\, z\bigg\}
&=\frac{1}{\sqrt{z}}\,\mathscr{K}_{0}\big\{x^{1/2}\,f(x)\, ;\, z\big\}\\
\label{r2:t1:2}
\mathscr{L}_{2}
\bigg\{\frac{1}{y^{2}}\,\mathscr{L}_{2}\Big\{f(x)\, ;\, \frac{1}{2y}\Big\}\, ;\, z\bigg\}
&=\mathscr{F}_{C}\Big\{\mathscr{G}\big\{x\,f(x)\, ;\, u\big\}\, ;\, z\Big\}
\intertext{and}
\label{r2:t1:3}
\mathscr{K}_{0}\big\{x^{1/2}\,f(x)\, ;\, z\big\}&=
\sqrt{z}\,\mathscr{F}_{C}\Big\{\mathscr{G}\big\{x\,f(x)\, ;\, u\big\}\, ;\, z\Big\} \cdot
\end{align}
\end{rem}
\begin{rem}\label{r3:t1} 
If we let $\nu=1/2$ in our Corollary \ref{c2:t1} and than use the formula \eqref{hankel:sin} and the definition \eqref{fsin} of the Fourier sine transform, we obtain
\begin{align}
	\label{r3:t1:1}
	\mathscr{L}_{2}
	\bigg\{\mathscr{L}_{2}\Big\{f(x)\, ;\, \frac{1}{2y}\Big\}\, ;\, z\bigg\}
	&=\frac{1}{2z^{3/2}}\,\mathscr{K}_{1}\big\{x^{3/2}\,f(x)\, ;\, z\big\}\\
	\label{r3:t1:2}
	\mathscr{L}_{2}
	\bigg\{\mathscr{L}_{2}\Big\{f(x)\, ;\, \frac{1}{2y}\Big\}\, ;\, z\bigg\}
	&=\frac{1}{2z}\,\mathscr{F}_{S}\Big\{u\,\mathscr{G}\big\{x\,f(x)\, ;\, u\big\}\, ;\, z\Big\}
	\intertext{and}
	\label{r3:t1:3}
	\mathscr{K}_{1}\big\{x^{3/2}\,f(x)\, ;\, z\big\}&=
	\sqrt{z}\,\mathscr{F}_{S}\Big\{u\,\mathscr{G}\big\{x\,f(x)\, ;\, u\big\}\, ;\, z\Big\}\cdot
\end{align}
\end{rem}
\begin{cor}\label{c3:t1} 
If the integrals involved converge absolutely, then we have
\begin{align}
\notag
\int_{0}^{\infty}\frac{1}{y}\,\exp\bigg(-\frac{z^{2}}{8y^{2}}\bigg)
\,{\rm I}_{\frac{\nu}{2}} \bigg(-\frac{z^{2}}{8y^{2}} \bigg)\,&\mathscr{L}_{2}\big\{f(x)\, ;\, y\big\}\,dy
\\
\label{c3:t1:1}
&=\int_{0}^{\infty} x\,f(x)\,{\rm I}_{\frac{\nu}{2}}\Big(\frac{z\,x}{2}\Big)\,
{\rm K}_{\frac{\nu}{2}}\Big(\frac{z\,x}{2}\Big)\,dx\\
\notag
\int_{0}^{\infty}\frac{1}{y}\,\exp\bigg(-\frac{z^{2}}{8y^{2}}\bigg)
\,{\rm I}_{\frac{\nu}{2}} \bigg(-\frac{z^{2}}{8y^{2}} \bigg)\,&\mathscr{L}_{2}\big\{f(x)\, ;\, y\big\}\,dy
\\
\label{c3:t1:2}
&=z^{-1/2}\,\mathscr{H}_{\nu}\,\Big\{u^{-1/2}\,\mathscr{G}\big\{x\,f(x)\, ;\, u\big\}\, ;\, z\Big\}
\\
\intertext{and}
\label{c3:t1:3}
\int_{0}^{\infty} x\,f(x)\,{\rm I}_{\frac{\nu}{2}}\Big(\frac{z\,x}{2}\Big)\,
{\rm K}_{\frac{\nu}{2}}\Big(\frac{z\,x}{2}\Big)\,dx&=z^{-1/2}\,\mathscr{H}_{\nu}\,\Big\{u^{-1/2}\,\mathscr{G}\big\{x\,f(x)\, ;\, u\big\}\, ;\, z\Big\}\cdot
\end{align}
\end{cor}
\begin{pf}
The proof of the Corollary \ref{c3:t1} is analogous to the previous Corollary \ref{c2:t1}. The assertions \eqref{c3:t1:1}, \eqref{c3:t1:2}, and \eqref{c3:t1:3} are obtained by putting
\begin{equation}
\label{c3:t1:p1}
g(u)=\frac{{\rm J}_\nu(z\,u)}{u}
\end{equation}
in our Theorem \ref{t1} and by using the known formulas \citet[p. 174, Entry (g)]{Gl} and \citet[p. 185, Entry (29)]{E1}. 
\qed
\end{pf}

The following corollary contains an identity involving $\mathscr{L}_{2}$-transform, the Glasser transform,  the $\mathscr{E}_{1}$-transform defined by
\begin{equation}
\label{e11}
\mathscr{E}_{1}\big\{f(x) \, ;\, y \big\}=\int_0^\infty \exp\big(x\,y\big)\,{\rm E}_1(x\,y\big)\,f(x)\,dx,
\end{equation}
introduced in \citet[p. 1377, Eq. (1.1)]{BDY:1}, $\mathscr{E}_{2,1}$-transform defined by
\begin{equation}
\label{e21}
\mathscr{E}_{2,1}\big\{f(x) \, ;\, y \big\}=\int_0^\infty x\,\exp\big(x^2\,y^2\big)\,{\rm E}_1(x^2\,y^2\big)\,f(x)\,dx,
\end{equation}
introduced in \citet[p. 1557, Eq. (1.1)]{BDY:2} and the Widder transform defined by 
\begin{equation}
\label{wid}
\mathscr{P}\big\{f(x) \, ;\, y \big\}=\int_0^\infty \frac{x\,f(x)}{x^2+y^2}\,dx
\end{equation}
introduced by Widder \citet{Wd}. The function ${\rm E}_{1}(x)$ is the second member of a family of functions defined by
\begin{equation}
\label{sn}
{\rm E}_{n}(x)=\int_{1}^{\infty} \frac{\exp(-x\,t)}{t^{n}}\,\quad n=0,1,\ldots\, .
\end{equation}
The functions $\text{E}_{n}(x)$ were introduced by  Schl\"omilch. The function ${\rm E}_{1}(x)$ in the definitions \eqref{e11} and \eqref{e21} of the $\mathscr{E}_{2,1}$ is related in a simple way to exponential integral function:
\begin{equation}
\label{s1:ei}
{\rm E}_{1}(x)=-\text{\rm Ei}(-x),
\end{equation}
where the exponential integral function is defined by
\begin{equation}
\label{ei}
\text{\rm Ei}(x)=\int_{-\infty}^{x} \frac{\exp(t)}{t}\,dt.
\end{equation}

\begin{cor}\label{c4:t1} 
If the integrals involved converge absolutely, then we have
\begin{equation}
\label{c4:t1:1}
\mathscr{E}_{2,1}\bigg\{\frac{1}{y}\,\mathscr{L}_{2}\big\{f(x)\, ;\, y\big\}\, ;\, z\bigg\}=
\sqrt{\pi}\,
\mathscr{P}\Big\{\mathscr{G}\big\{x\,f(x)\, ;\, u\big\}\, ;\, z\Big\}\cdot
\end{equation}
\end{cor}

\begin{pf}
The assertion \eqref{c4:t1:1} immediately follows upon putting
\begin{equation}
\label{c4:t1:p1}
g(u)=\frac{1}{u^{2}+z^{2}}
\end{equation}
the identity \eqref{t1:2} of our Theorem \ref{t1} and by using the known formula \citet[p. 185, Entry (29)]{E1}.
\end{pf}
\section{Illustrative Examples}
\numberwithin{equation}{section}

An interesting illustration for the identity \eqref{l1} asserted by Lemma \ref{l1} is contained in the following example. 
\begin{exmp}\label{e1}
Suppose that $|z|>|y|$. Then
\begin{gather}
\label{e1:1}
\mathscr{L}_{2}\bigg\{\frac{1}{u}\,\exp\big(z^{2}\,u^{2}\big)\,\text{E}_{1}\big(z^{2}\,u^{2}\big)\,  ;\,  y\bigg\}=\sqrt{\pi}\,\frac{\pi-2\,\arcsin(y/z)}{\sqrt{z^{2}-y^{2}}},
\\
\label{e1:2}
\mathscr{L}\bigg\{\frac{1}{\sqrt{u}}\,\exp\big(z\,u\big)\,\text{E}_{1}\big(z\,u\big)\, ;\, y\bigg\}=\sqrt{\pi}\,\frac{\pi-2\,\arcsin\big(\,\sqrt{y/z}\,\big)}{\sqrt{z-y}},\\
\label{e1:3}
\mathscr{E}_{2,1}\bigg\{\frac{1}{u}\,\exp\big(-y^{2}\,u^{2}\big)\, ;\, z\bigg\}=\sqrt{\pi}\,\frac{\pi-2\,\arcsin(y/z)}{\sqrt{z^{2}-y^{2}}},
\intertext{and}
\label{e1:4}
\mathscr{E}_{1}\bigg\{\frac{1}{\sqrt{u}}\,\exp\big(-y\,u\big)\, ;\, z\bigg\}=\sqrt{\pi}\,\frac{\pi-2\,\arcsin\big(\,\sqrt{y/z}\,\big)}{\sqrt{z-y}}\cdot
\end{gather}
\end{exmp}
\begin{pf}
We put
\begin{equation}
\label{e1:p1}
f(x)=\frac{1}{x^{2}+z^{2}}\cdot
\end{equation}
Using the known result \citet[p. 10, Entry (47)]{Ap} we find that
\begin{equation}
\label{e1:p2}
\mathscr{G}\Big\{\frac{x}{x^{2}+z^{2}}\, ;\, y\Big\}=\frac{\pi-2\,\arcsin(y/z)}{2\,\sqrt{z^{2}-y^{2}}}
\cdot
\end{equation}
Using the relationship \eqref{l2:l} between the Laplace transform and the $\mathscr{L}_{2}$-transform and the known formula \citet[p. 17, Entry (5)]{PBM}, we obtain
\begin{equation}
\label{e1:p3}
\mathscr{L}_{2}\Big\{\frac{1}{x^{2}+z^{2}}\, ;\, u\Big\}=\frac{1}{2}\,\mathscr{L}\Big\{\frac{1}{x+z^{2}}\, ;\, u^{2}\Big\}
=\frac{1}{2}\,\exp\big(z^{2}\,x^{2}\big)\,\text{E}_{1}\big(z^{2}\,x^{2}\big)\cdot
\end{equation}
Substituting the results \eqref{e1:p2} and \eqref{e1:p3} into the identity \eqref{l1:1} of our Lemma \ref{l1}, we obtain the asserted formula \eqref{e1:1}. 

From the relationship \eqref{l2:l} we deduce the assertion \eqref{e1:2}. The assertions \eqref{e1:3} and \eqref{e1:4} follow from the definitions \eqref{e11} and \eqref{e21} of the $\mathscr{E}_{1}$-transform and $\mathscr{E}_{2,1}$-transform, respectively.
\qed
\end{pf}
The following illustration involves the error function defined by
\begin{equation}
\label{err}
\text{\rm Erf}(x)=\frac{2}{\sqrt{\pi}}\,\int_{0}^{x} \exp\big(-x^{2}\big)\,dx,
\end{equation}
and the Dawson integral defined by
\begin{equation}
\label{daw}
\text{\rm daw}(x)=\int_{0}^{x} \exp\big(t^{2}-x^{2}\big)\,dt.
\end{equation}
The Dawson integral and the error function are related via the identity
\begin{equation}
\label{err:daw}
\text{\rm daw}(x)= \frac{-i\,\sqrt{\pi}}{2}\,\exp\big(-x^{2}\big)\,\text{\rm Erf}(i\,x)
\end{equation}
(cf. \citet[p. 405, Eq. 42:0:1)]{SO}).
\begin{exmp}\label{e2}
We show
\begin{gather}
\label{e2:1}
\mathscr{L}_{2}\bigg\{\frac{1}{u^{2}}\,\exp\bigg(-\frac{z^{2}}{4u^{2}}\bigg)\,\text{\rm Erf}\Big(i\,\frac{z}{2u}\Big)\, ;\, y\bigg\}=\frac{\pi\,i}{2}\,\big[\text{\rm I}_{0}(z\,y)-\mathbf{L}_{0}(z\,y)\big]
\intertext{and}
\label{e2:2}
\mathscr{L}_{2}\bigg\{\frac{1}{u^{2}}\,\text{\rm daw}\Big(\frac{z}{2u}\Big)\, ;\, y\bigg\}=\frac{\pi^{3/2}}{4}\,\big[\text{\rm I}_{0}(z\,y)-\mathbf{L}_{0}(z\,y)\big],
\end{gather}
where $\text{\rm I}_{0}(x)$ denotes the modified Bessel function of the first kind of order zero and $\text{\bf \rm L}_{0}(x)$ denotes the modified Struve function of order zero. 
\end{exmp}
\begin{pf}
We put
\begin{equation}
\label{e2:p1}
f(x)=\frac{\sin(z\,x)}{x}\cdot
\end{equation}
Using the relationship \eqref{l2:l} and the known formula \citet[p. 154, Entry (36)]{E1}, we have
\begin{align}
\notag
\mathscr{L}_{2}\big\{x^{-1}\,\sin(z\,x)\,;\,u\big\}&=
\frac{1}{2}\,\mathscr{L}\big\{x^{-1/2}\,\sin(z\,x^{1/2})\,;\,u^{2}\big\}\\
\label{e2:p2}
&=
-\frac{i\,\sqrt{\pi}}{2u}\,\exp\Big(-\frac{z^{2}}{4u^{2}}\Big)\,\text{\rm Erf}\Big(i\,\frac{z}{2u}\Big)\cdot
\end{align}
Multiplying both sides of \eqref{e2:p2} by $1/u$ and then applying the $\mathscr{L}_{2}$-transform, we find that
\begin{equation}
\label{e2:p3}
\mathscr{L}_{2}\Big\{\frac{1}{u}\,\mathscr{L}_{2}\big\{x^{-1}\,\sin(z\,x)\,;\,u\big\}\,;\,y\Big\}
=
-\frac{i\,\sqrt{\pi}}{2}\,\mathscr{L}_{2}\Big\{\frac{1}{u^{2}}\,\exp\Big(-\frac{z^{2}}{4u^{2}}\Big)\,\text{\rm Erf}\Big(i\,\frac{z}{2u}\Big)\,;\,y\Big\}\cdot
\end{equation}
From the known formula \citet[p. 174, Entry (a)]{Gl} we have
\begin{equation}
\label{e2:p4}
\mathscr{G}\big\{\sin(z\,x)\,;\,y\big\}=
\frac{\pi}{2}\,\big[\text{\rm I}_{0}(z\,y)-\mathbf{L}_{0}(z\,y)\Big]\cdot
\end{equation}
Substituting the formulas \eqref{e2:p3} and \eqref{e2:p4} into the identity \eqref{l1} of our Lemma \ref{l1}, we obtain the desired result \eqref{e2:1}.

From the relationship \eqref{err:daw}  and the formula \eqref{e2:1} we deduce, the assertion \eqref{e2:2}.
\qed
\end{pf}
\begin{rem}\label{r1:e2} 
Using the relationship \eqref{l2:l} and setting $z=2a^{1/2}$ we can restate the formulas \eqref{e2:1} and \eqref{e2:2} as 
\begin{gather}
\label{r1:e2:1}
\mathscr{L}\bigg\{\frac{1}{u}\,\exp\bigg(-\frac{a}{u}\bigg)\,\text{\rm Erf}\Big(i\,\sqrt{\frac{a}{u}}\,\Big)\, ;\, y\bigg\}
=\pi\,i\,\big[\text{\rm I}_{0}\big(2\,\sqrt{a\,y}\big)-\mathbf{L}_{0}\big(2\,\sqrt{a\,y}\big)\big]
\intertext{and}
\label{r1:e2:2}
\mathscr{L}\bigg\{\frac{1}{u}\,\text{\rm daw}\Big(\,\sqrt{\frac{a}{u}}\,\Big)\, ;\, y\bigg\}=\pi\,\big[\text{\rm I}_{0}\big(2\sqrt{a\,y}\big)-\mathbf{L}_{0}\big(2\sqrt{a\,y}\big)\big].
\end{gather}
\end{rem}
\begin{exmp}\label{e3} 
Suppose that $\Re(z)>0$ and $\max\{0,-2\,\Re(\nu)\}<\Re(\mu)<1$. Then
\begin{gather}
\label{e3:1}
\int_{0}^{\infty} y^{-\mu-1}\,\exp\bigg(-\dfrac{z^{2}}{2y^{2}\,}\bigg)\,{\rm I}_{\nu}\bigg(\,\frac{z^{2}}{2y^{2}}\,\bigg)\,dy
=\frac{\Gamma\Big(\dfrac{1}{2}-\dfrac{\mu}{2}\Big)\,\Gamma\Big(\nu+\dfrac{\mu}{2}\Big)}{2\sqrt{\pi}\,z^{\mu}\,\Gamma\Big(\nu-\dfrac{\mu}{2}+1\Big)}
\intertext{and}
\label{e3:2}
\int_{0}^{\infty} u^{\mu-1}\,\text{\rm I}_{\nu}(z\,y)\,\text{\rm K}_{\nu}(z\,y)\,dy
=
\frac{\Gamma\Big(\dfrac{\mu}{2}\Big)\,\Gamma\Big(\dfrac{1}{2}-\dfrac{\mu}{2}\Big)\,\Gamma\Big(\nu+\dfrac{\mu}{2}\Big)}{4\sqrt{\pi}\,z^{\mu}\,\Gamma\Big(\nu-\dfrac{\mu}{2}+1\Big)}
\end{gather}
(cf. \citet[p. 13, Eq. (4.8)]{KSY}).
\end{exmp}

\begin{pf}
Putting  
\begin{equation}
\label{e3:p1}
f(x)=\frac{\text{\rm J}_{2\nu}(2zx)}{x}
\end{equation}
into the identity \eqref{c1:t1:1} of Corollary \ref{c1:t1}, we obtain
\begin{equation}
\label{e3:p2}
\int_{0}^{\infty} y^{-\mu}\,\mathscr{L}_{2}\bigg\{\frac{\text{\rm J}_{2\nu}(2zx)}{x}\, ;\, y\bigg\}\,dy
=\frac{1}{2}\,\Gamma\Big(\frac{1}{2}-\frac{\mu}{2}\Big)
\,\int_{0}^{\infty}x^{\mu-1}\,\text{\rm J}_{2\nu}(2zx)\,dx,\end{equation}
Utilizing the formulas \eqref{c3:l1:p2} and \citet[p.326, Entry (1)]{E1} together with \eqref{e3:p2} we obtain the assertion \eqref{e3:1}.

Similarly, using the formula \eqref{c3:l1:1} and \citet[p.326, Entry (1)]{E1} together with \eqref{c1:t1:3} we deduce the second assertion \eqref{e3:2} of our Example \ref{e3}.

\end{pf}
We conclude this investigation by remarking that many other infinite
integrals can be evaluated in this manner by applying the above Lemma, the 
above Theorem, and their various corollaries and consequences considered 
here.




\end{document}